\documentclass{amsart} 
\usepackage{amssymb, amscd,url}
\numberwithin{equation}{section}

\def\CC{{\mathbb C}}
\def\DD{{\mathbb D}}

\def\HH{{\mathbb H}}

\def\PP{{\mathbb P}}
\def\QQ{{\mathbb Q}} 
\def\RR{{\mathbb R}} 
\def\VV{{\mathbb V}} 
\def\ZZ{{\mathbb Z}}

\def\half{\tfrac{1}{2}}
\def\G{{\Gamma}}
\def\g{{\gamma}}
\def\i {{\sqrt{-1}}}

\def\an{{\rm an}} 
\def\bb{{\rm bb}}

\def\reg{{\rm reg}}

\def\st{{\rm st}}
\def\ss{{\rm ss}}
\def\sing{{\rm sing}}

\def\bb{{bb}}

\def\bs{\backslash}
\def\bss{{\bs\!\bs}}
\def\pt{{\bullet}}
\def\eps{\epsilon}

\def\Fcal{{\mathcal F}} 
 
\def\Hcal{{\mathcal H}}

\def\Kcal{{\mathcal K}}
\def\Lcal{{\mathcal L}}

\def\Ocal{{\mathcal O}}

\def\Xcal{{\mathcal X}}

\newcommand\Coker{\operatorname{Coker}} 
 
\newcommand\Gr{\operatorname{Gr}}
\newcommand\Hom{\operatorname{Hom}}

\newcommand\im{\operatorname{Im}}
\newcommand\Ker{\operatorname{Ker}}

\newcommand\proj{\operatorname{Proj}}
\newcommand\res{\operatorname{Res}}

\newcommand\PO{\operatorname{PO}}

\newcommand\SL{\operatorname{SL}}

\newcommand\Orth{\operatorname{O}}

\newtheorem{theorem}{Theorem}[section]
\newtheorem{lemma}[theorem]{Lemma}
\newtheorem{proposition}[theorem]{Proposition}
\newtheorem{propdef}[theorem]{Proposition-Definition}

\newtheorem{corollary}[theorem]{Corollary}

\theoremstyle{definition}
\newtheorem{definition}[theorem]{Definition}

\theoremstyle{remark} 
\newtheorem{remark}[theorem]{Remark}

\begin{document}

\title{On period maps that are open embeddings}
\author{Eduard Looijenga}
\author{Rogier Swierstra}
\email{looijeng@math.uu.nl, swierstr@math.uu.nl}
\address{Mathematisch Instituut\\
Universiteit Utrecht\\
P.O.~Box 80.010, NL-3508 TA Utrecht\\
Nederland}
\thanks{Swierstra is supported by the Netherlands Organisation for Scientific
Research (NWO)}


\begin{abstract} 
For certain complex projective manifolds (such as K3 surfaces
and their higher dimensional analogues, the complex symplectic projective 
manifolds) the period map takes values in a locally symmetric variety of 
type IV. It is often an open embedding and in such cases it has been observed that 
the image is the complement of a locally symmetric divisor. We 
explain that phenomenon and get our 
hands on the complementary divisor in terms of geometric data. 
\end{abstract}

\maketitle

\section*{Introduction}
The period map assigns to a complex nonsingular projective variety  the isomorphism type of the polarized  Hodge structure on its primitive cohomology in some fixed degree; it therefore goes from a moduli space which parametrizes 
varieties to a moduli space which parametrizes polarized Hodge structures. 
The latter is always a locally homogeneous complex manifold which comes with an invariant metric. This manifold need not be locally symmetric, but 
when it is, then the period map is often a local isomorphism. 
Prime examples are (besides the somewhat tautological case of polarized 
abelian varieties) projective K3 surfaces and more generally, projective complex
symplectic manifolds:  the period map then takes values a locally symmetric variety of type IV or a locally symmetric subvariety thereof such as a ball quotient.
In many  of these cases, the period map can be proved to be even an open immersion
and we then find ourselves immediately wondering what the complement of the
image might be. It has been observed that this is almost always a 
locally symmetric  arrangement complement, that is, the complement of a finite union of locally symmetric hypersurfaces and this was the main reason for one of us to
develop a compactification technique (that generalizes the Baily-Borel theory) for such complements. 

In this paper we approach the issue from the other, geometric, side, by 
offering (among other things) an 
\emph{explanation} for the said observation. Our main result in this 
direction is Corollary
\ref{cor:main} which gives sufficient conditions  
for the image of a period map to be a locally symmetric arrangement 
complement once it is known to be an open embedding. The point is that 
in many cases of interest these conditions are known to be satisfied or are verified with 
relative ease. 
This is particularly so with K3 surfaces and we illustrate that with the 
Kulikov models and a case that was already analyzed by K\=ondo \cite
{kondo}, namely the quartic surfaces that arise as a cyclic cover 
of the projective plane along a quartic curve. 
But our chief motivation is to apply this to cases where the image
of the period map has not yet been established, such as cubic fourfolds 
(whose period map takes 
values in a locally symmetric variety of type IV of dimension 20), 
and more specifically, to those that arise as a cyclic cover of $\PP^4$ 
along a cubic hypersurface (whose period map takes values in a locally 
complex hyperbolic variety of dimension 10). This last case will be the 
subject of a subsequent paper. 

It is worth noting that our main result also provides a criterion for 
the surjectivity of the 
period map; this refers of course to the case when it implies that the locally symmetric divisor 
is empty. What we like about this criterion is that it is \emph{nonintrusive}:
it deals with degenerations as we encounter them in nature, while
leaving them untouched and, as Lemma \ref{lemma:typegeom} will testify,
does not ask us to subject the degenerating family to some kind of 
artificial processing (of which is unpredictable how and when it will 
terminate). 

The title of this a paper is explained by the fact that almost 
every significant result in this paper has as a hypothesis that some period map
is an open embedding. This can be hard to verify,
even if we know the period map to be a local isomorphism. Nevertheless the technique
developed here sometimes allows us to derive this hypothesis from
much weaker assumptions (see Remarks \ref{rem:weakinj} and \ref{rem:injex}).
\\

Finally a few details about the way this paper is organized. 
The first section is about polarized 
variations of Hodge structure `of type IV' defined on the generic point 
of an irreducible variety.
It begins with a simple Lemma
(namely \ref{lemma:limitperiod}), which contains one of the basic 
ideas on which this paper is built. We develop this idea with the help of 
our notion of \emph{boundary extension} of  a polarized variation of 
Hodge structure, leading up to the main results of this section, 
Propositions \ref{prop:supset}  and \ref{prop:supsetbis}. Such a 
boundary extension yields  rather precise information on the image 
of the classifying map of certain variations of polarized Hodge structure 
(also the injectivity of the classifynig map is addressed, here).
Section 2 brings the discussion 
to the geometric stage, specifically, to the setting of geometric 
invariant theory. It contains our main result
Corollary \ref{cor:main}. One can see this corollary at work in the next 
section, where it is being applied to (families of) K3 surfaces. 
We here also make the connection with the
classical theory of Kulikov models and we illustrate the  power
of the theory by applying it to the (well-understood) case of the moduli space  of 
quartic curves. In the appendix we show that the Baily-Borel 
type compactifactions of locally symmetric arrangement 
complements come with a natural boundary extension. In the appendix we show
that the natural compactification of a locally symmetric arrangement complement
comes with a natural boundary extension.
\\

In what follows a Hodge structure is always assumed to be defined over $\QQ$
and if the Hodge structure is polarized, then this also applies to the 
polarizing form. This allows us to identify a  polarized Hodge structure of weight zero 
with its dual. A variation of Hodge structure is not only supposed to have an underlying local system of $\QQ$-vector spaces,  but the latter is also assumed to contains some unspecified (flat) $\ZZ$-sublattice. 

If $\HH$ is a variation of Hodge structure over a complex manifold $M$, then
we denote by $\Fcal^\pt (\HH)$ the Hodge flag on 
$\Ocal^\an_M\otimes_{\CC_M} \HH$. We use the traditional superscript ${}^*$ as the generic way to indicate the dual of an object, where of course the meaning depends on the category 
(which we usually do not mention since it is clear from the context). So $\HH^*$ is the $\CC_M$-dual of $\HH$ as a variation of Hodge structure, but
$\Fcal^p (\HH)^*$ is the $\Ocal^\an_M$-dual of $\Fcal^p (\HH)$.

\section{A limit theorem for certain period maps}\label{sect:locdisc}

\subsection*{A local discussion} The following notion will be used throughout.

\begin{definition}
A Hodge structure of even weight $2k$ is said to be \emph{of type IV} if is
polarized and its Hodge number $h^{k+i,k-i}$ is  1 for $|i|=1$ and $0$ for $|i|>1$. 
\end{definition}

This implies that its quadratic form has signature $(h^{k,k},2)$ or 
$(2, h^{k,k})$ according to whether $k$ is even or odd. Notice that tensoring with
the Tate structure $\CC (k)$ turns a weight $2k$ Hodge structure of type IV into a weight $0$ Hodge structure 
of that type. 

If we are only given a vector space $H$ with $\QQ$-structure endowed 
with a $\QQ$-quadratic form of signature $(m,2)$, then to make it 
a weight zero Hodge structure $H$ of type IV is to choose a complex 
line $F^1\subset H$ with the property that $F^1$ is isotropic for the quadratic form and negative definite for the
hermitian form: if $\alpha\in  F^1$ is a generator, then  we want 
$\alpha\cdot\alpha=0$ and  $\alpha\cdot\bar \alpha<0$. The subset of $H$ 
defined by these two conditions has two connected components, interchanged 
by complex conjugation, one of which we will denote by $H_+$. The 
projectivization $\PP(H_+)$ is the symmetric space for the orthogonal group 
of $H(\RR)$ and is known as a symmetric domain of type IV. This explains our 
terminology, since  $\PP(H_+)$  classifies Hodge structures of that type.

We shall also consider situations where the group $\mu_l$ of $l$th roots of unity  with $l\ge 3$ acts on $H$ in such a manner that the eigenspace $H_\chi$ for the tautological character $\chi :\mu_l\subset\CC^\times$, has hyperbolic
signature $(m,1)$ (relative to the hermitian form).  We observe that $H_\chi$ is
isotropic for the quadratic form, for given $\alpha,\beta\in H_\chi$ and $\zeta\in\mu_l$ of order $l$,  then $\alpha\cdot\beta=\zeta\alpha\cdot\zeta\beta=\zeta^2\alpha\cdot\beta$
and hence $\alpha\cdot\beta=0$. The open subset $H_{+,\chi}$ of $H_\chi$ defined by $\alpha\cdot\bar\alpha<0$ has as its projectivization $\PP(H_{\chi,+})$ 
a complex ball in $\PP(H_\chi)$. This ball is the symmetric  space of the unitary group of the Hermitian form  on $H_\chi$  and is apparently also the classifying space for Hodge structures of type IV with a certain $\mu_l$-symmetry. Since $H_{\chi,+}$ is connected, it is contained in exactly one component $H_+$ as above.

We now state and prove a simple lemma  and discuss its consequences. 

\begin{lemma}\label{lemma:limitperiod}
Let be given a normal complex analytic variety $B$, a smooth Zariski open-dense 
subset $j: B^\circ\subset B$ and a variation of weight zero 
polarized Hodge structure $\HH$ of of type IV over $B^\circ$. 
Suppose further given a point $o\in B-B^\circ$, 
a subspace $V\subset H^0(B^\circ,\HH)$ (the space of flat sections
of $\HH$) and a generating section
$\alpha\in H^0(B^\circ, \Fcal^1(\HH))$  such that 
\begin{enumerate}
\item[(i)] for every $v\in V$, the function $s\in B^\circ\mapsto  v\cdot \alpha(s)$ extends holomorphically to $B$ and
\item[(ii)]$-\alpha (s)\cdot\bar\alpha (s)\to\infty$ as $s\in B^\circ\to o$.
\end{enumerate}
Then any limit of lines $F^1(s)$ for  $s\to o$ lies in $V^\perp$ in 
the following sense: if $U\subset B^\circ$ is simply connected with $o$ 
in its closure and  $p\in U$ is some base point, then 
\[
P: s\in U\mapsto [F^1(s)\subset \HH(s)\cong \HH(p)]\in \PP(\HH(p)),
\]
has the property that any accumulation point of $P(s)$ for $s\to o$ 
lies in the orthogonal complement of $V$ in $\HH (s)$.
\end{lemma}

\begin{proof} 
If $k:=\dim V$, then choose a basis $e_1,\dots ,e_N$  of $\HH|_U$ such
that $e_{k+1},\dots ,e_N$ are perpendicular to $V$.
So if we write $\alpha (s)=\sum_{i=1}^N 
f_i(s)e_i$, with each $f_i$ holomorphic on $U$, then the first 
assumption implies 
that $f_1,\dots ,f_k$  extend holomorphically across $o$. The second 
assumption tells us that $\sum_r |f_r|\to \infty$  as $s\to o$. This
implies that any accumulation point of $[f_(s):\cdots :f_N(s)]$ for $s\to o$
has its first $k$ coordinates zero.
\end{proof}

\begin{remark}
Condition (i) essentially comes down to requiring that the image of 
\[
V\subset (j_*\HH)_o\cong (j_*\HH^*)_o\subset (j_*(\Ocal^\an_{B,o}\otimes\HH^*))_o\to (j_*\Fcal^1(\HH)^*)_o
\]
be contained in an principal $\Ocal^\an_{B,o}$-submodule: $v\in V$ sends $\alpha$ to
$v\cdot\alpha\in \Ocal^\an_{B,o}$ and so the $\Ocal^\an_{B,o}$-module spanned by the image of $V$ in $(j_*\Fcal^1(\HH)^*)_o$ can be identified with the ideal spanned by such
functions in $\Ocal^\an_{B,o}$. Conversely, if the image of the above map is contained
in a principal $\Ocal^\an_{B,o}$-submodule of $(j_*\Fcal^1(\HH)^*)_o$, 
generated by $\alpha$, say, then $\alpha$ satisfies (i) on a neighborhood of $o$ in $B$.
Notice that if the image of $\cdot \alpha$ in $V^*$  (under evalution in $o$) 
is nonzero, then $V$ itself generates a pricipal $\Ocal^\an_{B,o}$-submodule of 
$(j_*\Fcal^1(\HH)^*)_o$. 

Condition (ii) says that on the given principal submodule (which is generated by $\alpha\mapsto 1$), the Hodge norm goes to zero at $o$ (so that on the norm on  its dual goes to infinity).
\end{remark}

\subsection*{Types of degeneration}
We return to the situation of Lemma \ref{lemma:limitperiod}. Let $p\in B^\circ$ be 
fixed base point and let us write $H$ for the vector space underlying $\HH(p)$
while retaining its $\QQ$-structure and the polarizing quadratic form. 
The latter has signature $(\dim H-2,2)$. We recall that the set of $\alpha\in H$ with 
$\alpha\cdot\alpha =0$,
$\alpha\cdot\bar\alpha <0$ has two connected components, one 
of which, denoted $H_+$,  contains a generator of $F^1(p)$. Its 
projectivization $\PP(H_+)$ is a symmetric domain for the orthogonal group 
of $H(\RR)$; it is also the domain for a classifying map of $\HH$, 
for we can think of $P$ as taking values in  $\PP(H_+)$. 
This is also a good occasion to recall that the boundary of $\PP(H_+)$
in $\PP(H)$ decomposes naturally into \emph{boundary components}:
a boundary component  is given by a nontrivial isotropic subspace $J\subset H$
defined over $\RR$ (so $\dim J=\{ 1,2\}$): the corresponding boundary component
is the $\PP(J)$-interior of $\PP(J)\cap \PP(H_+)^-$. So if $\dim J=1$, it is the singleton $\PP(J)$ and if $\dim J=2$ we get
an open half space on $\PP(J)$. The only incidence relations between these boundary components come from inclusions: if the closure of the boundary component
attached to $J$ meets the one associated to $J'$, then $J'\subset J$.

Assume now that $V$ is defined over $\RR$. We are  given that there exists a sequence $(\alpha_i)_i$ in $H_+$ with the lines $(\CC\alpha_i)_i$  
converging to some line $F_\infty\subset H$. According to Lemma \ref{lemma:limitperiod}
we have $F_\infty\perp V$ so that $[F_\infty]\in \PP (V^\perp)\cap \PP(H_+)^-$. The following is clear.

\begin{lemma}\label{lemma:types}
Let $V_0\subset V$ denote the nilspace of the quadratic form.
If the image of $V$ in $(j_*\HH)_o$ is defined over $\RR$,  
then  we are in one of the following three cases:
\begin{enumerate}
\item $V_0=0$. Then $V$ is positive definite, $V^\perp$ has signature $(\dim V^\perp-2,2)$  and $\PP(H_+)\cap \PP (V^\perp)$ is a nonempty
(totally geodesically embedded) symmetric subdomain of $\PP(H_+)$.
\item $\dim V_0=2$. Then $V$ is  positive semidefinite, $V^\perp$ is
negative semidefinite and $\PP(V^\perp)\cap \PP(H_+)^-=
\PP(V_0)\cap \PP(H_+)^-$.
\item $\dim V_0=1$ and  $\PP(V^\perp)\cap \PP(H_+)^-$ is the union of the boundary components of $\PP(H_+)$ that have $\PP(V_0)$ in their closure.
\end{enumerate}
\end{lemma}

So in the last two cases, $\PP(V^\perp)$ does not meet  
$\PP(H_+)$, but does meet its boundary. We would like to be able to
say that in case (1) $[F_\infty]\in \PP(H_+)$, that in the other two cases
$[F_\infty]$ lies in the boundary component defined by $V_0$ and that
$V$ is also positive semidefinite in case (3). We shall see that we come close
to fulfilling these wishes if we assume:
\begin{enumerate}
\item[(i)] the subspace $V\subset (j_*\HH)_o$ is defined over $\QQ$ and
\item[(ii)] the image of $\alpha$ in $V^*$  under evaluation in $o$ 
is nonzero so that it spans a line $F\subset V^*$.
\end{enumerate}
But then the discussion is no longer elementary, as we need to invoke the
mixed Hodge theory of one-parameter degenerations.
For this purpose we make a base change over the open unit disk $\Delta\subset \CC$,
which sends $0$ to $o$ and $\Delta^*$ to $B^\circ$. 
We assume here simply that $H$ is the fiber of a base point 
in the image of $\Delta^*$. Let $T:H\to H$ denote the monodromy operator of 
the family over $\Delta$. It is known that some positive power 
$T^k$ is unipotent. Since we can arrive at this situation by a 
finite base change, we assume that this is already the case. 

Now $T-1$ is nilpotent (we shall see that in the present case its 
third power is zero) and hence  $N:=\log T= -\sum_{k\ge 1}\frac{1}{k} 
(1-T)^k$ is a finite sum and nilpotent also. Notice that $N$ will be a rational element of 
the Lie algebra of the orthogonal group of $H$. If $N$ is not the 
zero map, then there exist linearly independent  $\QQ$-vectors 
$e,u$ in $H$ with $e\cdot e=e\cdot u=0$  such that 
\[
N(a)= (a\cdot e)u-(a\cdot u)e.
\]
Hence $T$ lies canonically in the natural one-parameter subgroup (take $w=1$)
\[
\exp (wN)(a)=a+w(a\cdot e)u-w(a\cdot u)e-\tfrac{1}{2}w^2(u\cdot u)(a\cdot
e)e,\quad w\in\CC.
\]
We have three cases:
\begin{enumerate}
\item[(I)] $N=0$,
\item[(II)] $N\not=0=N^2$, which means that $(u\cdot u)= 0$ and
\item[(III)] $N^2\not=0=N^3$, which means that $(u\cdot u)\not= 0$. 
\end{enumerate}
Let $J$ denote  the span of $e$ and $u$ and $J_0$ the nilspace
of $J$ (so $J_0$ equals $J$ resp.\ $\CC e$ in case II resp.\ case III). 
Notice that $V\subset\Ker (N)=J^\perp$ and that $V^*$ is a quotient of
$\Coker (N)=H/J$.

We kill the monodromy by counteracting it as follows. 
A universal cover $\widetilde{\Delta^*}\to\Delta^*$ of $\Delta^*$ 
can be taken to be the upper half plane with coordinate $w$  so that 
the covering projection is given by $s=\exp (2\pi\sqrt{-1}w)$ and 
$w\mapsto w+1$ generates the covering group. The variation of Hodge 
structure over $\Delta^*$ is given by a holomorphic map 
$P: \widetilde{\Delta^*}\to \PP(H_+)$
with $P(w+1)=TP(w)$. Then $\exp(-w\log N)P(w)$ only depends on
$s=\exp (2\pi\sqrt{-1}w)$ and so we get a holomorphic map 
$\phi : \Delta^*\to \PP(H)$. Schmid \cite{schmid}  proved that the 
latter extends holomorphically across $0\in\Delta$.
The line $F_{\lim}$ defines the Hodge filtration of a mixed Hodge structure 
$H_{\lim}$ on $H$  whose weight filtration is the Jacobson-Morozov 
filtration $W_\pt$ defined by $N$. This makes $N$ a morphism of 
mixed Hodge structures $N:H_{\lim} \to H_{\lim} (-1)$. 
The pure weight subquotients are polarized with the help of $N$. 
It follows that $H/J$ has a natural mixed Hodge structure with
$\pi_J(F_{\lim})$ defining the Hodge filtration. In particular
$\pi_J(F_{\lim})$ is nontrivial. Since $N$ is the zero map in $H/J$, it 
follows that if  take limits in  $H/J$ instead, we find that
\[
\lim_{\im (w)\to\infty}\pi_J(F_w)=\pi_J(F_{\lim}).
\]
where $\pi_{J}: H\to H/J$ is the projection. Notice that if the image 
of the lefthand side under the projection $H/J\to V^*$ is nonzero, then
it must equal $F$. So if in addition $V^*$ has a mixed Hodge structure for which 
$H/J\to V^*$ is a MHS-morphism, then its Hodge filtration is given by $F$.

Let us now go through the three cases.

(I)  $H_{\lim}$ is pure of weight zero and the period map 
factors through an analytic map $\Delta\to \PP(H_+)$ which takes in $0$ the 
value $[F_{\lim}]$. So in this case $V$ is positive definite and
is a Hodge substructure of $H_{\lim}$ which is pure of bidegree $(0,0)$.
In particular, the line $F\subset V^*$ has no Hodge theoretic significance. 

(II) Here $0=W_{-2}\subset W_{-1}=J \subset J^\perp=W_0 \subset W_1=H$
and $F_{\lim}$ projects nontrivially in  $H/J^\perp\cong J^*$; the latter
has a Hodge structure of weight 1 and $J^\perp/J$ which is pure of 
bidegree $(0,0)$. 
The line $F\subset V^*$ is the image of $F_{\lim}$ in $V^*$, but 
we cannot conclude that $V^*$ thus acquires a Hodge structure for which
$H/J\to V^*$ is a MHS-morphism unless we know that
$V_0=J$ (a priori, $0\not= V_0\subset J$); in that case
$(V/V_0)^*$ is pure of bidegree $(0,0)$ and $V_0^*$ is of weight $1$.

(III) Then $0=W_{-3}\subset W_{-2}=W_{-1}=J_0 \subset 
J_0^\perp=W_{0}=W_1 \subset W_2=H$ and $F_{\lim}$ projects isomorphically onto
$H/J_0^\perp\cong J_0^*$. The latter
is isomorphic to $\CC(-1)$ and polarized by 
the form $(N^2a\cdot b)$. Since we have $(N^2a\cdot a)=-(u\cdot u)(a\cdot e)^2$
it follows that $(u\cdot u)<0$. So $J$ is negative semidefinite
and $J^\perp$ is positive semidefinite (and hence so is $V$). In particular, $V_0=J_0$. Thus $V^*$ acquires a Hodge structure
with $(V/V_0)^*$ pure of bidegree $(0,0)$ and $V_0^*$ of bidegree $(1,1)$.

We sum up our findings and use the occasion to make a definition:

\begin{propdef}\label{propdef:lperiod}
Assume that in the situation of Lemma \ref{lemma:limitperiod}, $V\subset (j_*\HH)_o$ is defined over $\QQ$ and
that the image of $\alpha$ in $V^*$  under evaluation in $o$ 
spans a nonzero line $F\subset V^*$. Then $V$ is positive 
semidefinite and  $V^*$ has a  mixed Hodge structure characterized as follows:
$(V/V_0)^*$ is of bidegree $(0,0)$, $V_0^*$ has pure weight equal to its dimension
and if $V_0\not=0$, then $F=F^1(V^*)$. 
It has the property that if we make a base change over an analytic curve  $\Delta\to B$ whose  special point goes to $o$ and whose generic point to $B^\circ$, 
then the natural maps $V\to H_{\lim}$ and $H_{\lim}\to V^*$ 
are morphisms of mixed Hodge structure, unless
$\dim V_0=1$ and the base change is of type II.

 If this last case never occurs, we
say that $V\subset (j_*\HH)_o$ is a mixed \emph{Hodge} subspace.
\end{propdef}

It is obvious that $(j_*\HH)_o$ itself is a mixed Hodge subspace.

\subsection*{A global version} Suppose now that we are given a  
irreducible normal variety $S$ and  a variation of polarized 
Hodge structure $\HH$ of type IV weight zero over a Zariski open-dense subset 
$j:S^\circ\subset S$. Let $S^f\subset S$ denote the set of $s\in S$ where 
$j_*\HH$ has finite monodromy. This is a 
Zariski open subset which contains $S^\circ$. We put  
$S_\infty:=S-S^f\subset S$. We choose a base point $p\in S^\circ$ and let $H$ and $H_+$ have the same meaning as before. We have  a monodromy 
representation $\pi_1(S^\circ, p)\to \Orth (H)$, whose image is the \emph{monodromy group} $\Gamma$ of $\HH$. It preserves $H_+$ and defines an unramified  $\G$-covering $\widetilde{S}^\circ\to S^\circ$ on which is defined 
the classifying map $P: \widetilde{S}^\circ\to \PP(H_+)$.  This map  is $\G$-equivariant. The $\G$-covering $\widetilde{S}^\circ\to S^\circ$ extends canonically 
as a ramified  $\G$-covering $\widetilde{S}^f\to S^f$ with normal total 
space. We know that the classifying map then extends as a 
complex-analytic $\G$-equivariant map  $P: \widetilde{S}^f\to \PP(H_+)$.  

Recall that the Baily-Borel theory asserts among other things that in case
$\G$ is arithmetic,  $\G\bs \PP(H_+)$ admits a natural projective completion
$\G\bs \PP(\widehat{H}_+)$ whose boundary is Zariski closed (so that
$\G\bs \PP(H_+)$ is in a natural manner a 
quasiprojective variety). Under mild assumptions, we are in that situation: 

\begin{lemma}\label{lemma:arithm}
If $S$ is complete, $\HH$ has regular singularities  along  
$S_\infty$ and $P: \widetilde{S}^f\to \PP(H_+)$ is an open map, then $\G$ is arithmetic and $P$  descends to an open morphism $S^f\to \G\bs \PP(H_+)$ in the quasiprojective category.
\end{lemma}
\begin{proof}
To prove that $\G$ is arithmetic, choose an arithmetic  $\G'\supset\G$
(since $\G$ stabilizes a lattice, we can take for $\G'$ the orthogonal 
group of that lattice). Then $P$ determines an analytic map $S^f\to\G'\bs \PP(H_+)$. 
Domain and range have algebraic compactifications (namely $S$ and the Baily-Borel 
compactification 
$\G'\bs \PP(\widehat{H}_+)$) and since the singularities of this map have no essential 
singularities, its graph in $S\times\G'\bs\PP(\widehat{H}_+)$ has algebraic closure. 
This implies that $S^f\to\G'\bs \PP(H_+)$ has finite (positive) degree. Since that map factors through $\G'\bs \PP(H_+)$, it follows that $[\G :\G']$ 
is finite. This proves that $\G$ is arithmetic also and that $P$  descends to an 
open morphism  $S^f\to \G\bs \PP(H_+)$ in the quasiprojective category.
\end{proof}
          
Suppose now the variation of polarized Hodge structure $\HH$ comes 
with a (fiberwise) faithful action of $\mu_l$ so that $\HH_\chi$ has 
hyperbolic signature and hence the classifying map induces a morphism 
taking values in $\PP (H_{\chi,+})$. We then have the following 
counterpart of the above lemma (whose proof is similar to the case treated
and therefore omitted): 

\begin{proposition}\label{prop:supsetballbis}
If $S$ is complete, $\HH$ has regular singularities along $S_\infty$ and $P: \widetilde{S}^f\to \PP(H_{\chi,+})$ is an open map, then $\G$ is arithmetic and $P$  
descends to an open morphism $S^f\to \G\bs \PP(H_{\chi,+})$ in the quasiprojective 
category.
\end{proposition}

We ready now ready to introduce the notions that are central to this paper.

\begin{definition}\label{def:ext}
We call a constructible subsheaf $\VV$ of $ j_*\HH$ a \emph{boundary extension} 
for $(S,\HH)$ if it is defined over $\QQ$, its quotient sheaf has support 
on $S_\infty$ and the image of the sheaf map
\[
\VV\subset j_*\HH\cong j_*\HH^*\subset j_*(\Ocal_{S}\otimes \HH^*)\to j_*\Fcal^1(\HH)^*
\]
generates an invertible $\Ocal_S$-submodule (whose restriction to
$S^\circ$ will be $\Ocal_{S^\circ}\otimes\Fcal^1\HH)^*$) such that its 
(degenerating) norm on defined by the polarization of $\HH$ tends to zero 
along $S_\infty$. If in addition the stalks of $\VV$ define mixed Hodge subspaces
of the stalks of $j_*\HH$ (in the sense of Proposition-definition \ref{propdef:lperiod}) then we say that $\VV$ is the boundary extension of is of \emph{Hodge type}. If $j_*\HH$  is a boundary extension, then we say that  $\HH$ is \emph{tight} on $S$.
\end{definition}
 
Clearly, a tight boundary extension is of Hodge type.
We shall denote the $\Ocal_S$-dual of the image of $\VV$ in $ j_*\Fcal^1(\HH)^*$
by $\Fcal$. It is an extension of $\Fcal^1(\HH)$ cross $S$ as a line bundle and has the property that the norm on $\Fcal$ relative to the polarization of  $\HH$ tends to infinity along $S_\infty$. So this reproduces the situation of Proposition-definition \ref{propdef:lperiod} stalkwise.
We shall refer to $\Fcal$ as the \emph{Hodge bundle} of the boundary extension.

\begin{remark} 
It can be shown that boundary Hodge extension carries in a natural manner the structure 
of a polarized Hodge modules in the sense of M.~Saito \cite{saito}. 
In particular, if we are given a boundary extension $\VV$ of  $(S,\HH)$ 
as in Definition \ref{def:ext}, then we can stratify $S$ into smooth 
subvarieries such that the restriction of $\VV$ to every stratum is 
a variation of polarized mixed Hodge structure.
\end{remark}

Let $\VV\subset  j_*\HH$ be a boundary extension for $(S,\HH)$. 
According to Lemma 
\ref{lemma:limitperiod},  for every $s\in S-S^f$, $\VV_s$ determines a $\G$-orbit 
in the Grassmannian $\Gr (H)$ of $H$. We denote by $\Kcal_s\subset \Gr (H)$ the collection of orthogonal complements of these subspace (which is also a $\G$-orbit).
It has a type (1, 2 or 3) according to the distinction made in Lemma \ref{lemma:types}. If we stratify $S_\infty$ into connected strata in such a manner that $\VV$ is locally constant on each 
stratum, then $s\mapsto \Kcal_s$ is constant on strata and so 
$\Kcal :=\cup_{s\in S-S^f} \Kcal_s$ is a finite union of 
$\G$-orbits in $\Gr (H)$.
We decompose according to type:
\[
\Kcal =\Kcal_1\cup\Kcal_2\cup\Kcal_3.
\]
The members of $\Kcal_2\cup\Kcal_3$ do not meet $H_+$, whereas the collection of linear sections 
$\{ \PP(K_+)\}_{K\in \Kcal_1}$ is locally finite on $\PP(H^+)$, because $\G$ preserves a lattice in $H$ (see \cite{looijenga}). So if we put
\[
H_+^\circ :=H_+ -\cup_{K\in \Kcal_1}K_+,\quad 
\PP(H_+^\circ) :=\PP(H_+) -\cup_{K\in \Kcal_1}\PP(K_+),
\]
then $H_+^\circ$ is open in $H_+$ and $\PP(H_+^\circ)$ is open in  $\PP(H_+)$. 
It is known that for every linear subspace $K\in H$ of signature $(n,2)$, the  image
of $\PP(K_+)$ in $\G\bs\PP(H_+)$ is a closed subvariety (see \cite{looijenga} for a proof). So  $\G\bs \PP(H^\circ_+)$ is Zariski open in $\G\bs \PP(H_+)$. 

\begin{proposition}\label{prop:supset}
Suppose that $S$ is complete, $\HH$ has regular singularities along $S_\infty$
and $P:\widetilde S_f\to \PP(H_+)$ is an open map.
Let also be given a boundary extension of $\HH$ across $S$ with associated
collection of type 1 subspaces $\Kcal_1\subset\Gr (H)$.
Put  $H_+^\circ:=H_+-\cup_{K\in \Kcal_1}K_+$ and denote by  $H_+^\circ\subset H_+^\diamond\subset H_+$ the hyperplane arrangement
complement of the collection of \emph{hyperplanes} $K\in \Kcal_1$ for 
which the image of $P$ is disjoint with $\PP(K_+)$ (so that $P$ maps to 
$\PP(H_+^\diamond)$). 
Then  $P(\widetilde S_f)$ contains  $\PP(H_+^\circ)$ and 
$\PP(H_+^\diamond)-P(\widetilde S_f)$ is everywhere of codim
$\ge 2$. 

If moreover $S_\infty$ is everywhere of codim $\ge 2$ in $S$,
$\Fcal$ is ample and  $P$ is injective,  
then the $\CC$-algebra of meromorphic automorphic forms
$\oplus_{d\ge 0}H^0(\PP(H_+^\diamond),\Lcal^{\otimes d})^\G$ is
finitely generated, $S$ can be identified with its $\proj$ and
$H_+^\circ=P(\widetilde S^f)= H_+^\diamond$. 
\end{proposition}
\begin{proof}
Under these hypotheses Lemma \ref{lemma:arithm} applies so that
the monodromy group $\G$ is arithmetic and 
$\G\bs P: S^f\to \G\bs \PP(H_+)$ is an open  morphism of varieties.
This implies that $\G\bs P$ has  Zariski open image and that the image
of $P$ is open-dense in $\PP(H_+)$.

We first prove that $P(\widetilde{S}^f)\supset\PP(H_+^\circ)$.
Let $\alpha\in\PP(H^\circ_+)$. Since $\alpha$  is in the 
closure of the image of $P$, there exists a sequence
$(\tilde s_i\in\widetilde{S}^f)_i$  for which $(\alpha_i:=P(\tilde s_i))_i$ 
converges to $\alpha\in\PP(H^\circ_+)$. Pass to a subsequence so that the 
image sequence $(s_i\in S^f)_i$ converges to some $s\in S$.
We cannot have $s\in S_\infty$, for then we must have $\alpha\in\PP(K)$ for 
some $K\in\Kcal_s$. So $s\in S^f$. Choose 
$\tilde s\in \widetilde{S}^f$ over $s$. It is clear that there exist 
$\g_i\in\G$ such that $(\g_i\tilde s_i)_i$ converges to $\tilde s$.
So the sequences  $(\g_i\alpha_i=P(\g_i\tilde s_i))_i$ and $(\alpha_i)_i$
converge in  $\PP(H_+)$ to $P(\tilde s)$ and  $\alpha$ respectively.
From the the fact that $\G$ acts properly discontinuously on $\PP(H_+)$, it 
follows that  a subsequence of  $(\g_i)_i$ is stationary, say equal to $\g$, so that
$\alpha=P(\g^{-1}\tilde s)\in P(\widetilde{S}^f)$. 

If $K\in\Kcal_1$ is such that the image of $P$ meets $\PP (K_+)$, then
$P$ meets $\PP(K_+)$ in a nonempty open subset and hence the same is true
for the image of the induced open embedding 
$\G\bs P :\widetilde S^f\to \G\bs\PP(H_+)$ with regard to the image 
of $\PP (K_+)$. In this last case, it does not matter whether we 
use the Hausdorff topology or the Zariski topology. This implies that  
$\PP(H_+^\diamond)-P(\widetilde S^f)$ has codim $\ge 2$ in $\PP(H_+^\diamond)$
everywhere.

The line bundle $\Fcal$ is the pull-back of the automorphic 
line bundle $\Lcal$ on $\G\bs\PP(H_+)$.
So if $S_\infty$ is everywhere of codim $\ge 2$ in $S$, then  
for every integer $d\ge 0$ we have an injection
\[
H^0(\PP(H_+^\diamond),\Lcal^{\otimes d})^\G\subset
H^0(S^f,\Fcal^{\otimes d})\cong H^0(S,\Fcal^{\otimes d}).   
\]
If in addition $P$ injective and $\Fcal$ is ample, then 
the displayed injection is an isomorphism and 
$S=\proj\left(\oplus_{d\ge 0}H^0(S,\Fcal^{\otimes d})\right)$
is a projective completion of $\G\bs\PP(H_+^\diamond)$. Hence
$\G\bs P$ will map $S^f$ isomorphically onto 
$\G\bs\PP(H_+^\diamond)$ and $\PP(H_+^\circ)=\PP(H_+^\diamond)$.
\end{proof}

\begin{remark}\label{rem:weakinj}
The requirement that $P$ be injective is (under the assumption that $P$ is
open) the conjunction of two conditions: that $P$ be a local isomorphism 
and that $P$ be of degree one. While it is often not difficult to verify 
the former, the latter is in general much harder to settle. Sometimes 
the injectivity  can be established along the way
if all the assumptions of Proposition \ref{prop:supset} are known to be
fulfilled (including the ampleness of $\Fcal$ and the codimension
condition for $S_\infty$) except that instead of $P$ being injective,  
we only know it to be a local isomorphism. To explain what we mean,
let us regard $\G\bs P$ as a rational map between two projective 
varieties: from $S$ to some compactification of $\G\bs \PP(H^\diamond_+)$ 
(for instance the Baily Borel compactification of $\G\bs \PP(H_+)$ or the
one we discuss below). If we are 
lucky enough to find a boundary point in the latter compactification 
such that the rational map is regular over a Hausdorff (or formal) 
neighborhood  of that point and is there an isomorphism, then
we may conclude $P$ that is injective so that
the last assertion of Proposition \ref{prop:supset} still holds.
\end{remark}

A natural completion of a variety of the form
$\G\bs H_+^\diamond$ was constructed in the two papers \cite{looijenga}.
This completion $\G\bs\PP(\widehat H_+^\circ)$ has the property that the automorphic line bundle on $\G\bs H_+^\diamond$ extends over it as an ample line bundle (in the orbifold sense, of course), so that the completion is in fact
projective. We recall its construction in the appendix and also show that it
comes with a natural boundary extension of the 
tautological variation of type IV Hodge structure over $\G\bs H_+^\diamond$. 
Following Theorem  \ref{thm:coext} it is tight over this compactification in case the boundary of the completion is everywhere of codim $\ge 2$.
The question comes up whether we get in the situation of the previous proposition
the same projective completion with the same boundary extension. Presumably the answer is always yes. Here we prove this to be so under some additional hypotheses.

\begin{proposition}\label{prop:supsetbis}
Assume we are in the situation of Proposition \ref{prop:supset}: 
$S$ is complete, $\HH$ has regular singularities on $S$, $\Fcal$ is ample on $S$ and $S_\infty$ everywhere of codim $\ge 2$ in $S$ and $P$  an open embedding.
If in addition $\dim S\ge 3$ and  the hermitian form on $H$  takes a positive value
on every two-dimensional intersection of
hyperplanes in $\Kcal_1$, then the isomorphism 
$S^f\cong \G\bs \PP(H_+^\circ)$ (asserted by Proposition
\ref{prop:supset}) extends to an isomorphism of $S$ onto the natural compactification 
$\G\bs\PP(\widehat H_+^\circ)$ of $\G\bs\PP(H_+^\circ)$. 
This isomorphism underlies an isomorphism of polarized variations of 
Hodge structure with boundary extension: the boundary extension of $\HH$ across $S$ is tight over $S$. 
\end{proposition}

\begin{proof}
The additional assumption and the fact that $\dim H=2+\dim \PP(H_+)=2+\dim S\ge 5$
imply that the boundary of the compactification of
$\G\bs\PP(H_+^\circ)\subset\G\bs\PP(\widehat H_+^\circ)$ is of codim $\ge 2$
everywhere. Since the automorphic line bundle $\Lcal$ on $\G\bs\PP(H_+^\circ)$ extends as an ample line bundle over the completion 
$\G\bs\PP(H_+^\circ)\subset\G\bs\PP(\widehat H_+^\circ)$, it follows that
\[
H^0(\PP(H_+^\circ),\Lcal^{\otimes d})^\G=
H^0(\PP(\widehat H_+^\circ),\Lcal^{\otimes d})^\G 
\]
and so 
\[
\G\bs\PP(\widehat H_+^\circ)=
\proj\left(\oplus_{d\ge 0}H^0(\PP(H_+^\circ),\Lcal^{\otimes d})^
\G\right).
\]  
This shows that $P$ induces an isomorphism $S\cong \G\bs\PP(\widehat
H_+^\circ)$. The assertion concerning the boundary extension follows 
from Theorem \ref{thm:coext}.
\end{proof}

If $\HH$ and its boundary extension over $S$ come with a 
(fiberwise) faithful action of $\mu_l$ so that $\HH_\chi$ has hyperbolic signature. 
Then we have the following counterpart of Propositions \ref{prop:supset} and 
\ref{prop:supsetbis} (proofs are omitted since they are similar to the case
treated): 

\begin{proposition}\label{prop:supsetball}
Suppose $S$ is complete, $\HH$ has regular singularities along $S_\infty$ and  
$P:\widetilde S^f\to \PP(H_{\chi,+})$ is an open map. 
Let $H_{\chi,+}^\circ\subset H_{\chi,+}^\diamond\subset H_{\chi,+}$ denote the arrangement complement of the collection of \emph{hyperplanes} $K_\chi$, $K\in \Kcal_1$, for which the image of $P$ is disjoint with $\PP(K_{\chi,+})$ (so that $P$ maps to $\PP(H_{\chi,+}^\diamond)$). Then $P(\widetilde S_f)$ contains $\PP(H_{\chi,+}^\circ)$ and has a complement in $\PP(H_{\chi,+}^\diamond)$ everywhere of codim $\ge 2$.

If moreover $P$ is injective, $S_\infty$ is everywhere of codim $\ge 2$ in $S$ and $\Fcal$ is ample, then in fact $H_{\chi,+}^\circ=P(\widetilde S^f)= H_{\chi,+}^\diamond$.

If in addition the hermitian form takes a positive value on every two-dimensional intersection of the hyperplanes of the form $K_\chi$, $K\in\Kcal_1$,  then the isomorphism $S^f\cong \G\bs \PP(H_{\chi,+}^\circ)$ 
extends to an isomorphism of $S$ onto the natural 
compactification $\G\bs\PP(\widehat H_{\chi,+}^\circ)$ of 
$\G\bs\PP(H_{\chi,+}^\circ)$ which underlies an isomorphism of polarized variations of 
Hodge structure with boundary extension.
In particular, the boundary extension is the direct image of 
the dual of $\HH$ on $S$. 
\end{proposition}

\section{The geometric context}\label{sect:hypo}

We begin with a definition:

\begin{definition}\label{def:hodgext}
Let $f:\Xcal\to S$ be a proper family of pure relative dimension $m$  
which is smooth over an open-dense connected $S^\circ\subset S$ and let 
$0<k\le m$ be an integer such that the fibers over $S^\circ$ have their 
cohomology in degree $2k$ of type IV. Denoting  by 
$f_\reg: (\Xcal/S)_\reg \to S$ the restriction of $f$ to the part where $f$ 
is smooth, then a \emph{geometric Hodge bundle for $f$  in degree $2k$} is  
a line bundle $\Fcal$ over $S$ together with an embedding of quasi-coherent 
sheaves $u:\Fcal\to (R^{2k}f_\reg^*f_\reg^{-1}\Ocal_S)(k)$ such that the induced
fiber maps $u(s): \Fcal (s)\to H^{2k}(X_{s,\reg},\CC)(k)$
enjoy the following two properties: 
\begin{enumerate}
\item[(i)] for every $s\in S$, $u(s): \Fcal (s)\to H^{2k}(X_{s,\reg},\CC)(k)$ is injective and 
\item[(ii)] when $s\in S^\circ$, the image of $u(s)$  is $H^{k+1,k-1}(X_s,\CC)(k)$ (so that $u$ identifies $\Fcal |S^\circ $ with 
$R^{k-1}f_*\Omega^{k+1}_{\Xcal /S}|S^\circ$).
\end{enumerate}
If $f$ is projective (so that the primitive part $\HH$ of $R^{2k} f_*f^{-1}\Ocal_S(k)|S^\circ$ is a polarized variation of Hodge structure), then we say that such a geometric 
Hodge bundle has \emph{a proper norm} if
\begin{enumerate}
\item[(iii)] the Hodge norm stays bounded on $\Fcal$ at $s$, precisely 
when $R^{2k}f_*\QQ_{\Xcal}|S^\circ$ has finite monodromy at $s$.
\end{enumerate}
\end{definition}

The justification of this definition is that a  geometric Hodge bundle 
$\Fcal$ with proper norm determines a boundary extension 
$\VV_f\subset j_*\HH$ as follows.
For every $s\in S$, there exists a neighborhood
$B$ of $s$ in $S$ and a retraction $r: X_B\to X_s$ which is
$C^\infty$-trivial over the smooth part $X_{s,\reg}$ of $X_s$ so that we get a 
$B$-embedding 
\[
\iota^s: X_{s,\reg}\times B\to (X_B/B)_\reg
\]
in the $C^\infty$-category. This embedding is unique up to $B$-isotopy.
In particular, we have a well-defined map 
$\iota^s_{s' *}:H_{2k}(X_{s,\reg},\CC)\to H_{2k}(X_{s',\reg},\CC)$ 
for every $s'\in B$. This defines a constructible
sheaf $\tilde \VV_f$ on $S$ whose stalk at $s$ is $H_{2k}(X_{s,\reg},\CC)(-k)$.  
(By means of fiberwise Poincar\'e duality (or rather Verdier duality), this sheaf can also be identified with $R^{2(m-k)}(f_\reg)_{ !}\CC (m-k)$.) The homomorphism $\HH^*\subset j^*\tilde\VV_f$ dualizes to a sheaf homomorphism
$\tilde\VV_f\to j_*\HH^*\cong j_*\HH$ and we let $\VV_f$ be the image of the latter.
So the stalk $\VV_{f,s}$ is the image of 
$H_{2k}(X_{s,\reg},\CC)(-k)$ in $\HH(s')^*\cong \HH (s')$, where $s'\in S^\circ$ is close to $s$.  
In order that $\VV_f$ defines a boundary extension, we need that the 
image of 
\[
\tilde\VV_f\to j_*\Hom_{\Ocal_{S^\circ}}
(j^*R^kf_*\Omega^{k}_{X/S}(k),\Ocal_{S^\circ})
\]
generates a line bundle. This is the case: that line bundle can be identified with the dual of $\Fcal$. The assumption that $\Fcal$ has proper norm ensures that
$\VV_f$ defines a boundary extension.\\

We now prepare for a geometric counterpart of the discussion in Section \ref{sect:locdisc}.
If $\HH$ is a variation of a polarized Hodge structure $\HH$ over a 
quasi-projective base, then according to Deligne \cite{deligne:hodge2}, $\HH$ 
is semisimple as a local system 
and its isotypical decomposition is one of  variation of a polarized Hodge 
structures.
In particular, the part invariant under monodromy is a polarized Hodge 
substructure and  so is its orthogonal complement. We refer to the latter 
as the \emph{transcendental part} of $\HH$.

Part of the preceding will be summed up by Corollary \ref{cor:main} below.
We state it in such a manner that it includes the ball quotient case, but
our formulation is dictated by the applications that we have in mind, rather than by
any desire to optimize for generality.

Let be given  
\begin{enumerate}
\item[(a)] a projective family $f:\Xcal\to S$  of pure  relative dimension $m$  with smooth general fiber and normal projective base, 
\item[(b)] an action on $f$ of a group of the form $G\times\mu_l$ where 
$G$ is (algebraic and) semisimple and $\mu_l$ is the group of  $l$th roots of unity ($l=1,2,\dots$), 
\item[(c)]  a $G\times\mu_l$-equivariant ample line bundle $\Fcal$ over $S$
(so that are defined the stable locus $S^\st$ and the semistable locus $S^\ss$ relative to the $G$-action on $\Fcal$) 
\item[(d)]  and an integer $0<k\le m$
\end{enumerate}
with the following properties:
\begin{enumerate}
\item[(i)] the  transcendental part $H$ of  a general fiber of 
$R^{2k}f_*\CC_\Xcal (k)$ is of type IV,
\item[(ii)] $\mu_l$ preserves the fibers of $f$ and acts on the generic stalk of
$R^{k-1}f_*\Omega^{k+1}_{\Xcal/S} (k)$ with tautological character $\chi:\mu_l\subset\CC^\times$,
\item[(iii)] if $H_\chi\subset H$ denotes the eigenspace of the tautological character $\chi:\mu_l\subset\CC^\times$, then the open subset $S^f\subset S$ 
where the local system $\HH_\chi$ defined by $H_\chi$ has finite monodromy
is contained in the stable locus  $S^\st$,
\item[(iv)] $G\bss (S^\ss-S^f)$ has codim $\ge 2$ everywhere in $G\bss S^\ss$,
\item[(v)] $\Fcal |S^\ss$ can be given the structure of a geometric Hodge bundle  for  $f_{S^\ss}$ in degree $2k$ such that it has proper norm,
\end{enumerate}
Then a period map $\widetilde S^f \to \PP(H_{\chi,+})$ is defined on the 
monodromy covering $\widetilde S^f$ of $S^f$. Since $G$ is semisimple, the $G$ action on $S^f$ lifts to an action of $\widetilde S^f$, perhaps after passing to a connected covering $\tilde G$ of $G$ such that the period map factors through 
$\tilde G\bs \widetilde S^f$. We finally require:
\begin{enumerate}
\item[(vi)] (Torelli property)  the map 
$\tilde G\bs \widetilde S^f\to  \PP(H_{\chi,+})$ through which the period map  factors is an open embedding.
\end{enumerate}

\begin{corollary}\label{cor:main}
Under these assumptions the monodromy group $\G$ of $\HH_\chi$ is arithmetic in its unitary resp.\ orthogonal group.  The collection hyperplanes of $\HH_\chi$
that appear as the kernel of  a member of the natural $\G$-orbit of maps 
$H_\chi\to H^{2k}(X_{s,\reg},\CC)_\chi$, where $s\in S^\ss-S^f$ has
a closed orbit in $S^\ss$, and which meet $H_{\chi,+}$ (that is, are of type I)
make up a $\G$-arrangement whose associated arrangement complement 
$\PP(H_{\chi,+}^\circ)\subset \PP(H_{\chi,+})$ is the exact image  of the period map so that there results an isomorphism $G\bs S^f\cong \G\bs \PP(H_{\chi,+}^\circ)$. If $\dim S\ge 3$ (or $\ge 2$ in case $l\ge 3$), then the associated
boundary extension of $\HH$ over $S^\ss$ is tight.
\end{corollary}
\begin{proof}
It is known that a variation of polarized Hodge structure of geometric origin has regular singularities. Since $G\bss S^\ss=\proj\left(\oplus_{d\ge 0} H^0\left( S,\Lcal^{\otimes d}\right)\right)^\G$ is a projective variety which carries the $G$-quotient of $\Lcal$ as an ample line bundle in the orbifold sense, the hypotheses of Propositions \ref{prop:supset} and \ref{prop:supsetball} are fulfilled after this passage to $G$-quotients and the corollary follows.
\end{proof}

The conditions (i)-(iv) and (vi)  of Corollary \ref{cor:main} are often easily checked 
in practice or are known to hold---it is usually the verification of (v) 
that requires work. This motivates the following definition.

\begin{definition}\label{def:bdyvar}
Let $X\subset \PP^N$ be a projective variety of pure dimension $m$ 
and $k$ a positive integer $\le m$ such that $X$ admits a projective
smoothing in $\PP^N$ whose general fiber has cohomology in degree $2k$ 
of type IV.   
We shall call a one dimensional subspace $F\subset H^{2k}(X_\reg,\CC)$ 
a \emph{residual Hodge line} if every such smoothing admits a geometric 
Hodge bundle in degree $2k$ which has $F$ as closed fiber
(this $F$ will  be unique, of course).  If in addition, the geometric 
Hodge bundle has proper norm, then we say that $X$ is a 
\emph{boundary variety in dimension $2k$}  and call  $(X,F)$ a \emph{boundary pair 
in degree $2k$}. 
\end{definition}

A smoothing of such a boundary variety  will have a type. Mixed Hodge theory allows us to read off this type from the pair $(X,F)$ without any reference to a smoothing:

\begin{lemma}\label{lemma:typegeom}
Let $(X,F)$ be a boundary pair  in degree $2k$ and let $w(F)$
be such that $F$ embeds in the weight $2k+w(F)$ subquotient
of $H^{2k}(X_\reg,\CC)$. Then $w(F)\in\{ 0,1,2\}$ and the type of any smoothing
is $w(F)+1$.
\end{lemma}
\begin{proof}
If $f:(\Xcal, X)\to (\Delta, 0)$ is a projective smoothing of any projective variety 
$X$ over the unit disk, then we have a natural map 
$H^0(\Delta^*,R^{\pt}f_*\QQ_\Xcal)\to H^{\pt}(X_\reg,\QQ)$, and hence also 
a map $H^0(\widetilde{\Delta^*},R^{\pt}f_*\QQ_\Xcal)\to H^{\pt}(X_\reg,\QQ)$.
It is known that if  we give  $H^0(\widetilde{\Delta^*},R^{\pt}f_*\QQ_\Xcal)$ the limiting mixed Hodge structure and $H^{\pt}(X_\reg,\QQ)$ the usual mixed Hodge structure, then this becomes a MHS-morphism. This fact, combined with the discussion  following 
Lemma \ref{lemma:types}  leading up to the three cases I, II and III  immediately yields the assertion.
\end{proof}

\section{Examples}

If in the setting of Corollary \ref{cor:main} the general fiber is a K3-surface 
(and $k=1$), then a logical choice for $\Fcal$ is $f_*\omega_{\Xcal/S}$ and 
(v) then boils down to the statement that for $s\in S^\ss$, a generator 
of $H^0(X_s,\omega_{X_s})$ defines a nonzero class in $H^2(X_{s,\reg},\CC)$ 
and does not lift as a regular form to a resolution of $X_s$. However, this 
last condition tends  to be always fulfilled:

\begin{proposition}\label{prop:k3}
Let $X$ be a projective surface with $H^1(X,\Ocal_X)=0$ and with 
trivial dualizing sheaf $\omega_X$. Then every smoothing of $X$ is a K3 surface.
A generator of $\omega_X$ is square integrable on $X$ 
precisely when $X$ is a K3 surface  with only rational double 
point singularities (so that $X$ is smoothable with finite 
local monodromy group). 
\end{proposition}
\begin{proof}
A smoothing of $X$ will be a surface with vanishing irregularity and trivial dualizing sheaf
and hence is a K3 surface.

Assume now that a generator of $\omega_X$ is square integrable on $X$. 
We first show that $X$ has isolated singularities. If 
not, then a general hyperplane section has a Gorenstein curve singularity 
with the property that a generator of its dualizing module is regular on the 
normalization. But this is impossible, since
the quotient of the two modules has dimension equal to  the Serre invariant, 
which vanishes only in the absence of a singularity.
Since $X$ has only isolated singularities, then these must be normal in view
of the fact that $X$ has no embedded components and $H^1(X,\Ocal_X)=0$. 
According to Laufer \cite{laufer}, the square integrability of $\alpha$ is then
equivalent to $X$ having only rational double 
points as singularities.  Then $X$ has a minimal resolution $\tilde X$ with
$H^1(\tilde X,\Ocal_{\tilde X})=0$ and trivial dualizing sheaf. 
This implies that
$\tilde X$ is a K3 surface.
\end{proof}

If $X$ is a smoothable projective surface as in the previous proposition 
(so with vanishing irregularity and trivial dualizing sheaf) and a generator
of $\omega_X$ is \emph{not} square integrable and has \emph{nonzero} image 
in $H^2(X_\reg,\CC)$, then it is clear from the definitions that $X$ is a boundary surface.
With the help of Lemma \ref{lemma:typegeom} we can specify the type. Let us first do this for some special cases.

\begin{proposition}\label{prop:ellcusp} 
A smoothable projective surface with vanishing 
irregularity and trivial dualizing sheaf which has a simple-elliptic resp.\
cusp singularity is a boundary surface of Hodge type 2 resp.\ 3.
\end{proposition} 
\begin{proof}
Let $X$ be such a surface and let $\alpha$ be a generating section of $\omega_X$
and let $X$ have at $p$ a  simple-elliptic or a cusp singularity.

If $p$ is simple elliptic, then the exceptional curve $C$ of the minimal resolution of the germ $X_p$ is a smooth genus one curve $C$ and  $\alpha$ has a simple pole of order one along $C$ with residue a generator of the dualizing sheaf of $C$.
It follows  that we have an embedding of Hodge structures $H^1(C,\QQ)(-1)\subset H^2(X_\reg,\QQ)$ whose image subsists in the coholomology of a limiting mixed Hodge structure of a smoothing. This implies that $X$ is of type 2.

If $p$ is a cusp singularity, then the exceptional curve $C$ of the minimal resolution of the germ $X_p$ is a is a cycle of rational curves and it is still true that $\alpha$ has a simple pole of order one along $C$ with residue a generator of the dualizing sheaf of 
$C$. The latter will have a nonzero residue at a singular point of $C$.
It follows  that we have an embedding  $(\QQ)(-2)\subset H^2(X_\reg,\QQ)$ whose image subsists in the coholomology of a limiting mixed Hodge structure of a smoothing. 
Hence  $X$ is of type 3.
\end{proof}

By essentially the  same argument we find:

\begin{proposition}\label{prop:k3case}
Let  $X$ be a smoothable projective surface with with vanishing irregularity 
and trivial dualizing sheaf. Suppose that for some resolution 
$\tilde X$ of $X$ there exists a connected component $C$ of $\tilde X-X_\reg$ 
that is smooth of genus one or a cycle of rational curves with only rational
double points, and has the property that a generator of $\omega_X$ has a
pole of exact order one along $C$. Then $X$ is a boundary surface of Hodge type
type 2 resp.\ 3 if $C$ is smooth resp.\ a rational cyle.
\end{proposition}

The customary approach to studying the period map of a degenerating 
family of K3 surfaces over the unit disk is to pass to a standard 
situation, called a \emph{Kulikov model}. This involves a combination of a 
base change and a blowing up over the central 
fiber (whose effect on the period map is merely a base change). The model in question has then the property that it has a smooth total space, a trivial dualizing sheaf
and a reduced normal crossing divisor as central fiber. But this process is highly nonunique and arriving at such a model tends to be nontrivial task. The interest of our approach lies in the fact that we leave such degenerations untouched and deal with them in the way they come: neither a central modification, nor usually  a base change are needed to extract useful information about  the limiting behavior of the period map. Let us nevertheless see how our method  deals with the Kulikov models.
These models fall into three classes according to the type of the central normal crossing surface $X$:
\begin{enumerate}
\item[(i)] $X$ is a smooth K3 surface.
\item[(ii)] $X$ is a chain of (at least two) smooth surfaces (the dual graph is an interval) whose double curves are mutually isomorphic  smooth genus one curves.  
The surfaces at the end are rational  and the surfaces in  between are elliptic ruled surfaces.
\item[(iii)] the normalization of $X$ consists of smooth rational surfaces and  the incidence complex of these surfaces is a triangulated two-sphere. The double curve meets every component in an anticanonical cycle.
\end{enumerate}
In case (i) there is nothing to do and Proposition \ref{prop:k3case} is applicable 
to the other two cases: a resolution of $X$ is given by normalization. In case (ii)
a member at the end of a chain in a rational surface $Y$ which meets the rest 
of the chain along a smooth genus one curve $C$ and  we have 
$\Ocal_Y\cong\omega_{X}\otimes \Ocal_Y=\omega_Y(C)$, so that $X$ is a boundary surface of Hodge type 2.
In case (iii), let $Y$ be any connected component of the normalization of $X$. 
Then $Y$ is rational and the preimage $C\subset Y$ of $X_{\sing}$ is a 
cycle of rational curves. So it is a boundary surface and according to
Proposition \ref{prop:k3case}, $X$ is a boundary surface of Hodge type 3.

\subsection*{A worked example: the moduli space of quartic curves}
The result we are going to discuss is essentially due to K\=ondo (\cite{kondo}, see also
\cite{looijenga2})--the point of what follows is merely to show how 
effectively it is reproduced by our method.
Let $G\subset\SL (4,\CC)$ be the stabilizer of the decomposition 
$\CC^4=\CC^3\oplus\CC$. This group is clearly isogenous to $\SL(3,\CC)\times\CC^\times$. It acts on the space of quartic forms in four complex variables 
$(Z_0, Z_1,Z_2,Z_3)$ of the shape $F(Z_0, Z_1,Z_2,Z_3)=f(Z_0, Z_1,Z_2)+\lambda Z_3^4$. Let us denote the projectivization of that space by $S$. This projective space supports a quartic surface $\Xcal\subset\PP^3_S$ with $\mu_4$-action. It is clear
that the open subset defined by $f\not= 0\not= \lambda$  parametrizes quartic surfaces in $\PP^3$ that are $\mu_4$-covers of quartic curves in $\PP^2$.
The $G$-semistable surfaces all lie in this open subset. The surface 
$X_s\subset\PP^3$ is 
$G$-(semi)stable if and only if the corresponding quartic plane curve $C_s\subset \PP^2$ is so relative to the quotient $\SL(3,\CC)$ of $G$. Following Mumford, $s\in S^\st$  precisely when  $C_s$ has only ordinary double points or (ordinary) cusps
(this means that $X_s$ has singularities of type $A_3$ or $E_6$) and 
$s \in S^\ss-S^\st$ has a closed orbit in $S^\ss-S^\st$ precisely when 
$C_s$ is a union of two conics $C', C''$ with $C'$ smooth,
$C''$ not a double line and for which either $C'$ and $C''$ meet in two 
points of multiplicity two, or  $C'=C''$. So if $C'\not=C''$, then either $C''$ is nonsingular and $C'$ and $C''$ have common tangents where they meet or 
$C''$ is the union of two distinct tangents of $C'$. These are represented by the
one-parameter family $(Z_1Z_2-Z_0^2)(Z_1Z_2-tZ_0^2)$, with $t\in\CC$
(but notice that we get the same orbit in $S$ for $t$ and $t^{-1}$). This
shows that $S^f=S^\st$. 

If $C=C_s$ is smooth, then $X=X_s$  is a polarized K3 surface of degree $4$ with $\mu_4$-action.  According to \cite{kondo} the  eigenspace $H^2(X,\CC)_\chi$ is hyperbolic of dim 7 and defined over $\QQ(\sqrt{-1})$. 

There is  a $G$-equivariant isomorphism of 
$f_*\omega_{\Xcal/S}|X^\ss\cong \Ocal_{S^\ss}(1)$ defined as follows: to every
$F\in \CC[Z_0, Z_1,Z_2,Z_3]$ as above defining a $G$-semistable surface $X$, we
associate the generator 
\[
\alpha_F:=\res_{X}\res_{\PP^3}\frac{dZ_0\wedge dZ_1\wedge dZ_2\wedge dZ_3}{F}.
\]
of $\omega_X$. Since  $F$ appears here with degree $-1$, we get the asserted 
line bundle isomorphism. It is evidently $G$-invariant.  
So if we take $\Fcal=\Ocal_S(1)$, then have all the data in place  for the verification of the properties of Corollary \ref{cor:main}.

\begin{proposition}
All the conditions of Corollary \ref{cor:main} are satisfied for  $\Fcal=\Ocal_S(1)$ and hence so is its conclusion: we get an isomorphism  $G\bss S^\ss\cong 
\G\bs\PP (\widehat H_{\chi,+}^\circ)$. Here $H_\chi$  is a vector space of dim 7 endowed with a hermitian form of hyperbolic signature defined over $\QQ(\sqrt{-1})$ and
$H_{\chi,+}$ is the arrangement complement  and the $\G$-arrangement 
defining $H_{\chi,+}$ consists of a single $\G$-orbit of hyperplanes 
of $H_\chi$ (these are perpendicular to a sublattice of $H_\ZZ$ of type 
$2I(-2)$ on which $\mu_4$ acts faithfully). We have $S^f=S^\st$ and $S^\ss-S^f$
has two strata: one stratum parametrizes unions of distinct conics, 
one of which is nonsingular and having exactly two points in common 
and  is of Hodge type 2, whereas the other is the orbit double nonsingular 
conics and is of type 1.  
\end{proposition}
\begin{proof}
Conditions (i) and (ii) of Corollary \ref{cor:main} are clearly satisfied.
We have already seen that $S^f=S^\st$, so that also
(iii) holds. We found that the orbit space 
$G\bss (S^\ss-S^f)$ is of dimension one, and hence (iv) is satisfied. 
Property (vi) is a consequence of  the Torelli theorem for K3-surfaces. 
So it remains to verify (v). We first check that $\Fcal$ indeed defines a 
geometric Hodge bundle over $S^\ss-S^f$.  A semicontinuity argument shows 
that it is enough to do this in the closed orbits of $G$. We verify 
this property and the remaining  (iv) at a point of such an orbit at the same time. 
We therefore assume that $C=C_s$ is the union of two conics $C'$ 
and $C''$ as described above.\\

\emph{The case  $C'=C''$.}
So $C$ is a nonsingular conic with multiplicity $2$.
Since $X=X_s\subset \PP^3$ is a $\mu_4$-cover over $\PP^2$ totally ramified 
along $C\subset \PP^2$, it consists of two copies $S',S''$ 
of $\PP^1\times \PP^1$ joined along their diagonal 
$D\subset \PP^1\times \PP^1$. A generating section $\alpha$
of $\omega_X$ is obtained as follows. If $z$ is the affine 
coordinate on $\PP^1$, then $\zeta:=(z'-z'')^{-2}dz'\wedge dz''$ 
extends as a regular $2$-form  on $\PP^1\times \PP^1-D$ with a 
pole of order $2$ along $D$. Then let $\alpha$ be on $S'$ equal to this 
form and on $S''$ minus this form.  

It is clear that $\zeta$ (and hence $\alpha$) is not square integrable.
We next verify that its cohomology class in $H^2(\PP^1\times \PP^1-D,\CC)$ 
is nonzero (so that the cohomology class of $\alpha$ in $H^2(X_\reg,\CC)$ 
is nonzero, also.) This we verify by evaluating the integral of $\zeta$
over the integral generator of $H_2(\PP^1\times \PP^1-D,\CC)$. A generator 
can be represented as follows: the algebraic cycle 
$[\PP^1]\otimes 1 -1\otimes[\PP^1]$ on $\PP^1\times\PP^1$ is homologous 
to a (nonalgebraic) cycle which avoids $D$. The homology takes place on a small neighborhood of
$(0,0)$: we let $\G_\eps$ be the sum of the holed Riemann spheres 
$(\PP^1-\Delta_\eps)\times \{0\}$ and $\{ 0\}\times (\PP^1-\Delta_\eps)$  
with opposite orientation (where $\Delta_\eps$ is the open $\eps$-ball 
centered at $0\in\PP^1$) with and the tube $T_\eps$ of
$(z',z'')\in\CC\times\CC\subset \PP^1\times\PP^1$ with $z'-z''=e^{\i\theta}$ 
and $z'+z''=te^{\i\theta}$, $-1\le t\le 1$. Then up to sign, we have
\begin{multline*}
\int_{\G_\eps}\zeta=\int_{T_\eps} \frac{dz'\wedge dz''}{(z'-z'')^2}=\\
=\int_{T_\eps} \frac{\half d(e^{\i\theta})\wedge d(te^{\i\theta})}{e^{2\i\theta}}=
 \int_{-1}^1\int_0^{2\pi} \half \i td\phi\wedge dt   =2\pi\i.
\end{multline*}
Notice that the selfintersection number of $\G_\eps$ is 
that of $[\PP^1]\otimes 1 -1\otimes[\PP^1]$ and hence equal to $-2$. 
Let  $\g$ be the corresponding class in $H_2(\PP^1\times \PP^1-D)$ lifted to $H_2(X_\reg)$. If $g\in\mu_4$ is a generator, then $g^2\g=-\g$ and $H_2(X_\reg)$ 
is freely generated by $\g$ and $g\g$. Thus, $H_2(X_\reg)$ is a $\mu_4$-module
isomorphic to the Gau\ss\  lattice $\ZZ[\sqrt{-1}]$ endowed with the quadratic form $-2\| z\|^2$. This lattice has no even overlattices and so  a copy of it gets 
\emph{primitively} and $\mu_4$-equivariantly embedded in the primitive homology of a nearby smooth quartic surface of the above type. This lattice is negative definite
and hence defines a hyperplane section of $\PP(H_{\chi,+})$; we
are in the type 1 case (compare \cite{kondo}).\\

\emph{The case $C'\not=C''$.}
Then the associated  $\mu_4$-cover $X$ has simple-elliptic singularities
of degree 2 at the two points of $C'\cap C''$. We can therefore invoke 
Proposition \ref{prop:ellcusp} to conclude that property (v) holds. 
Let us nevertheless do this in some detail. 
The two conics generate a pencil of which  the generic member completely decomposes in $X$. If we resolve the simple-elliptic singularities minimally, we get a surface $\tilde X$ with $\mu_4$-action that is obtained as follows.
Let $E$ be a smooth genus one curve with $\mu_4$-action isomorphic to $\CC/\ZZ[\sqrt{-1}]$ with its obvious $\mu_4$-action. This action has two distinct fixed points which we denote by $p',p''$. Now let $\mu_4$ act on $\PP^1$ via $\CC^\times\subset\PP^1$ and consider $\PP^1\times E$ with the diagonal $\mu_4$-action. 
Blow up $(0,p')$ and $(0,p'')$ in $\PP^1\times E$ and then blow down the strict transforms of $\PP^1\times\{p'\}$ and $\PP^1\times\{p''\}$. The result is a is a 
$\PP^1$-bundle $\bar X\to E$ with $\mu_4$-action. The fixed point set of the action is the union of the fibers over $p'$ and $p''$. The `zero section' has self-intersection $-2$ and the `section at infinity' has  self-intersection $+2$. Now choose $q\in E$.

If $q$ is not fixed by the $\mu_4$-action, then we let $\tilde X\to \bar X$ be the blow up the $\mu_4$-orbit of $(\infty,q)$. The sections $E_0, E_\infty$ 
of $\tilde X\to E$ at zero and infinity have both self-intersection $-2$ so that their contraction $\tilde X\to X$
yields a surface $X$ with two simple-elliptic singularities. The latter inherits a 
$\mu_4$-action whose orbit space is isomorphic to $\PP^2$. The fiber over $p'$ and $p''$ map to $C'$ and $C''$ and the fiber over $q$ maps to the union of the common tangents of $C'$ and $C''$. We have that $X_\reg\cong\tilde X-E_0-E_\infty$ and a straightforward calculation shows that  we have an exact sequence
\[
0\to \CC(-1)\to H^2(X_\reg,\CC)_\chi\to H^1(E_0,\CC)(-1)\to\ 0  ,
\]
in particular, $\dim H^2(X_\reg,\CC)_\chi=2$.
If $q$ is a fixed point, say $q=p''$, then we let $\tilde X\to \bar X$ be obtained by blowing up four times over $q$ (always on the strict transform of $\{\infty\}\times E$). This produces over $q$ a string of three $(-2)$-curves that can be contracted to $A_3$-singularity. The morphism $\tilde X\to X$ contracts also
the sections $E_0$ and $E_\infty$. Then 
$H^2(X_\reg,\CC)_\chi=H^2(\tilde X-E_0-E_\infty,\CC)_\chi\cong H^1(E_0,\CC)(-1)$, hence has dim 1.
In either case, a generating section of $\omega_X$ is obtained as follows:
if $\alpha_E$ is a nowhere zero regular differential $\alpha_E$ on $E$, 
then $z^{-1}dz\otimes \alpha_E$ defines a rational differential on $\tilde X$ 
with divisor $-E_0-E_\infty$ which is the pull-back of a unique generating 
section $\alpha$ of $\omega_X$. By Proposition
\ref{prop:ellcusp} we are in the type 2 case. 
\end{proof}

\begin{remark}\label{rem:injex}
The use of the Torelli theorem for K3 surfaces to establish property (v)
of Corollary \ref{cor:main} can be avoided if we use an argue as in Remark \ref{rem:weakinj} by taking
as compactification of $\G\bs \PP(H_+^\circ)$ the natural one: over 
the point of a one-dimensional stratum of the latter we find the 
simple-elliptic locus. The period map is near that locus is sufficiently 
well understood \cite{looij:se} for proving that along this locus we have a local isomorphism 
for the Hausdorff topology.   
\end{remark}

\section{Appendix: the boundary extension attached to an arrangement compactification}
We show here that the arrangement compactifications introduced in \cite{looijenga}
come with a natural boundary extension. The setting in which this construction is 
carried out is essentially the one that we are dealing with, for it needs 
the following input:
\begin{enumerate}
\item[(i)] a complex vector space $H$ with $\QQ$-structure of dimension 
$n+2\ge 4$ endowed with a nondegenerate quadratic form, also defined over 
$\QQ$ and of signature $(n,2)$, and a connected component $H_+$ of the set
of $\alpha\in H$ with $\alpha\cdot\alpha=0$ and $\alpha\cdot\bar \alpha<0$,
\item[(ii)] an arithmetic subgroup $\G$ of $\Orth (H)$ which stabilizes $H_+$ and 
\item[(iii)] a $\G$-invariant collection $\Hcal$ of hyperplanes of $H$, all defined 
over $\QQ$ and of signature $(n-1,2)$ and consisting of finitely many $\G$-equivalence classes. 
\end{enumerate}
We shall refer to these data as a \emph{$\G$-arrangement in $H_+$}  and we call
\[
H_+^\circ:=H_+ -\cup_{K\in\Hcal} K_+ \quad\text{resp.}\quad  
\PP(H_+^\circ):=\PP(H_+) -\cup_{K\in\Hcal} \PP(K_+).
\]
a \emph{$\G$-arrangement complement} in $H_+$ resp.\ $\PP (H_+)$. Given our
$\G$-arrangement $\Hcal$, we define an index set $\tilde\Sigma$ parametrizing  such subspaces that is a union 
\[
\tilde\Sigma:=\Kcal_1\cup \Kcal_{2}\cup(\cup_J\Sigma_J),
\]
of  which the first two members are self-indexing in the sense that these are 
collections of subspaces of $H$. 

We let $\Kcal_1$ be the collection of subspaces of $H$ with positive 
definite orthogonal complement that arise as an intersection of members 
of $\Hcal$ (this was denoted $\PO(\Hcal |\PP(H_+)$ in \cite{looijenga}).
If  $J\subset H$ is an isotropic plane defined over $\QQ$, then denote by  $K_J$ the common intersection of $J^\perp$ and all the members of 
$\Hcal$ which contain $J$ (so $K_J=J^\perp$ if no member of $\Hcal$ contains $J$).
This will be our collection $\Kcal_{2}$; it is clearly in bijective correpondence with the set isotropic planes defined over $\QQ$. 

Let now $J\subset H$ be an isotropic line defined over $\QQ$. Then 
$J^\perp/J$ is nondegenerate defined over $\QQ$ and has hyperbolic 
signature $(1,n)$. Let $e\in J(\QQ)$ be a
generator. The map $p_e: H-J^\perp\to (J^\perp/J)(\RR)$, 
$\alpha\mapsto \im\left((\alpha\cdot e)^{-1}\alpha\right)$ has
the property that $H_+$ is the preimage of a component $C_e$ of 
the quadratic cone in $J^\perp (\RR)/J(\RR)$ defined by $y\cdot y<0$.
It is more intrinsic to consider instead 
$p_J: H-J^\perp\to (J\otimes J^\perp /J)(\RR)$ given by
$\alpha\mapsto e\otimes\im\left((\alpha\cdot e)^{-1}\alpha\right)$.
So in the latter space there is a quadratic cone $C_J$ such that
$H_+=p_J^{-1}(C_J)$. (This map factors through $\PP (H_+)$ and gives rise
to the latter's realization as a tube domain of the first kind.)
Any member of $\Hcal$ which
contains $J$ defines hyperplane section of $C_J\cong C_e$. 
These hyperplane sections are locally finite on $C_J$ and decompose $C_J$ 
into a locally polyhedral cones. Denote this collection of cones by 
$\Sigma_J$. For any $\sigma\in \Sigma_J$, denote by
$K_\sigma$ the subspace of $J^\perp$ defined by the complex-linear span of 
$\sigma$ in $J^\perp/J$. So if no member of $\Hcal$ contains $J$, then 
$C_J$ is the unique member of $\Sigma_J$ and we then have $K_{C_J}=J^\perp$. 
This is in general not an injectively indexed collection of subspaces, for
it often happens that for distinct $\sigma,\sigma'$ we have 
$K_\sigma=K_{\sigma'}$. 
Notice that  $\Sigma_J$ is in a natural manner a partially ordered set.
The  disjoint union $\Kcal$ of the 
$\Sigma_J$, where $J$  runs over all the $\QQ$-isotropic lines in $H$,
$\Kcal_1$ and $\Kcal_{2}$ is then also partially ordered by taking 
the inclusion relation  between corresponding subspaces they define, 
except that on each $\Sigma_J$ we replace this by the inclusion relation 
between cones.
 
For every linear subspace $K\subset H$ defined over $\QQ$ for which
$\PP(K)$ meets the closure of $\PP(H_+)$ (so that $K$ is of type 1, 2 or
3), we denote by $\pi_K$ the projection $H\to H/K$ and by 
$\PP\pi_K: \PP (H)-\PP(K)\to \PP (H/K)$  its projectivization.
Consider the disjoint union
\[
\widehat H_+^\circ:=H_+^\circ \sqcup
\coprod_{\sigma\in\tilde\Sigma} \pi_{K_\sigma}(H_+^\circ)
\]
and its projectivization
\[
\PP(\widehat H_+^\circ):=\PP(H_+^\circ)\sqcup \coprod_{\sigma\in\tilde\Sigma} \PP\pi_{K_\sigma}\PP(H_+^\circ).
\]
In \cite{looijenga}-II, pp.\ 570 the latter is endowed with a 
$\G$-invariant Hausdorff topology which induces the given topology 
on the parts. The former
is not formally introduced there, but its definition is completely analogous for 
it is such that the obvious map  $\widehat H_+^\circ\to
\PP(\widehat H_+^\circ)$ is the formation of a $\CC^\times$-orbit space.
Both have the property that  the partial ordering prescribes the incidence 
relations. Notice that each of these spaces comes with a $\G$-invariant 
\emph{structure} sheaf of continuous complex valued functions,
namely the functions that are holomorphic on every member of the partition.
The $\G$-orbit spaces $\PP(\G\bs \widehat H_+^\circ)$ 
and $\G\bs \widehat H_+^\circ $ are Hausdorff and compact resp.\ locally 
compact and their structure sheaves make them normal analytic spaces 
for which the natural projection  
$\G\bs \widehat H_+^\circ\to\G\bs \PP(\widehat H_+^\circ )$ is a 
$\CC^\times$-bundle in the orbifold sense. The associated
orbifold line bundle (which we denote by $\Fcal$) turns out to be ample, so that 
$\G\bs \PP(\widehat H_+^\circ)$ is projective. 

We will need the  following lemma, which almost captures the topology of
$\widehat H_+^\circ$. 

\begin{lemma}\label{lemma:convergence}
If a sequence $(\alpha_i\in H_+^\circ)_i$ converges to $\alpha_\infty\in \pi_{K}(H_+^\circ)$, 
with $K$ a linear subspace indexed by $\tilde\Sigma$, then $\lim_{i\to \infty} \pi_{K}(\alpha_i)=\alpha_\infty$ and
$\lim_{i\to\infty} \alpha_i\cdot\bar\alpha_i=-\infty$.
\end{lemma}

This result is somewhat hidden in \cite{looijenga}-II, which makes it a bit 
hard to explicate.  The construction of  $\PP(\widehat H_+^\circ )$ involves 
an extension of $\PP(H_+)$ (which we will here denote by $\PP(\widehat H_+)$, 
but which is there denoted $\DD^{\Sigma(\Hcal)}$). We first notice that the linear subspaces which do not meet $H_+$ are indexed by $\Sigma:= \Kcal_{2}\cup(\cup_J\Sigma_J)$, so that 
we can form
\[
\widehat H_+=H_+\sqcup \coprod_{\sigma\in\Sigma}\pi_{K_\sigma}(H_+)\quad\text{and}\quad
\quad \PP(\widehat H_+)=\PP(H_+)\sqcup \coprod_{\sigma\in\Sigma}
\PP\pi_{K_\sigma}\PP(H_+).
\]
The topologies are described in \cite{looijenga}-II, p.\ 566. In case 
$\Hcal$ is empty, then $\widehat H_+$ is the Baily-Borel extension $H_+^\bb$ 
of $H_+$: the orbit space $\G\bs\PP (H^\bb_+)$ is the Baily-Borel compactification 
of $\G\bs\PP(H_+)$ and the $\CC^\times$-bundle $H^\bb_+\to \PP(H_+^\bb)$ is 
associated to the basic automorphic line bundle: if $\Lcal$ denotes the line 
bundle over $PP(H_+^\bb)$ obtained as a quotient of $\CC\times H^\bb_+$ by the $\CC^\times$ action
defined by $\lambda(z,\alpha)=(z\lambda^{-1},\lambda \alpha)$, then 
the $\G$-automorphic forms of degree $l$ are the continuous $\G$-equivariant 
sections of $\Lcal^{\otimes l}$ that are holomorphic on strata. 

There is an evident $\G$-equivariant map $\widehat H_+\to H_+^\bb$. 
It is such that the resulting map $\G\bs\PP(\widehat H_+)\to\G\bs\PP(H_+^\bb)$ is a morphism and involves a modification of the Baily-Borel boundary. This modification has the property that the
closure of the image of every $\PP(K_+)$, $K\in\Hcal$, is a $\QQ$-Cartier 
divisor. In other words, it can locally be given by a single equation.
Notice that the pull-back of the of the basic automorphic line bundle is the 
line bundle associated to the $\CC^\times$-bundle
$\widehat H_+\to \PP(\widehat H_+)$.

The relation between $\PP(\widehat H_+)$ and $\PP(\widehat H_+^\circ)$
is as follows: let $\widetilde{\PP} (H_+)$ be the blow up of
$\PP(\widehat H_+)$ obtained by blowing up the closures of the linear sections 
$\PP(K_+)\subset\PP(H_+)$, $K\in\Kcal_1$, in the order of increasing  dimension (where we of course take their strict transforms). Then there is a blowdown
$\widetilde{\PP} (H_+)\to \PP(\widehat H^\circ_+)$ obtained as follows: 
if $K\in\Kcal_1$ is of codimension $r$, then the exceptional 
divisor in $\widetilde{\PP} (H_+)$ associated to $K$ has a product structure with one factor equal to a modified $\PP^r$. We successively blow down onto these factors, starting with the factors of lowest dimension first with indeed 
$\PP(\widehat H^\circ_+ )$ as the final result. 
Let us denote by $\hat\Lcal$ the line bundle over $\PP(\widehat H^\circ_+)$ that is associated to $\widehat H_+\to \PP(\widehat H_+)$ (in the same way as $H^\bb_+\to \PP(H_+^\bb)$ is to $\Lcal$). On the common blowup $\widetilde{\PP} (H_+)$, 
the pull-back of $\hat\Lcal$ is the pull-back of $\Lcal(\Hcal):=\Lcal(\sum_{K\in\Hcal})$. This helps us to settle part of the proof of Lemma \ref{lemma:convergence}.

\begin{proof}[Proof of Lemma \ref{lemma:convergence}] 
Suppose that the $K$ appearing in the lemma is a member of $\Kcal_1$. We
regard $H_+^\circ$ as an open subset of the total space of the pull-back of
$\Lcal(\Hcal)$ to $\widetilde{\PP} (H_+)$. The assumption is that 
the sequence $\alpha_i$ converges in this total space to an element $\alpha_\infty$
that does not lie on the zero section. Notice that the sequence $([\alpha_i]\in\PP(H_+))_i$ converges to some $[\alpha]\in\PP (K_+)$. 
This means that if the sequence $\alpha_i$ is regarded as one which lies in the total space of $\Lcal$, then there exist scalars $\lambda_i\in\CC$ with $\lambda_i\alpha_i$ converging in $H_+$ to an element $\alpha\in H_+$ and 
$\lim_{i\to\infty}\lambda_i =0$. So $\lim_{i\to \infty}|\lambda_i|^2\alpha_i\cdot\bar\alpha_i=\alpha\cdot\bar\alpha< 0$ and hence $\lim_{i\to \infty}
\alpha_i\cdot\bar\alpha_i=-\infty$. 

It remaining cases (when $K$ is degenerate) involve a different type of argument
and are settled in Lemma \ref{lemma:convergence2} below.
\end{proof}

\begin{lemma}\label{lemma:convergence2}
If a sequence $(\alpha_i\in H_+)_i$ converges $\alpha_\infty\in \pi_{K}(H_+)$, 
with $K$ a linear subspace indexed by $\Sigma$, then 
$\lim_{i\to \infty} \pi_{K}(\alpha_i)=\alpha_\infty$ and 
$\lim_{i\to\infty} \alpha_i\cdot\bar\alpha_i=-\infty$.
\end{lemma}
\begin{proof}
Suppose first that $K=K_J$ for some isotropic plane $J$. 
Then $J(\RR)$ has a natural orientation characterized by the following 
property: if $(e_0,e_1)$ is an oriented basis for $J(\RR)$, then for all
$\alpha\in H^+$, we have $\im ((\alpha\cdot e_0)(\bar\alpha\cdot e_1))>0$.
Given such a basis $(e_0,e_1)$, we define a one parameter group
of orthogonal transformations as follows: for $\tau\in\CC$,
\[
\psi_\tau :H\to H,\quad  \psi_\tau(\alpha)=\alpha +\tau\left((\alpha\cdot e_0)e_1
-(\alpha\cdot e_1)e_0\right).
\]
Notice that $\psi_\tau$ leaves $H/J$ (and hence $H/K$) pointwise fixed.
A small computation shows that $\psi_\tau$ preserves the quadratic form 
on $\HH$ and that
\[
\psi_\tau(\alpha)\cdot\overline{\psi_\tau(\alpha)}=
\alpha\cdot\bar\alpha -4\im(\tau) \im \left((\alpha\cdot e_0)(\bar\alpha\cdot
e_1)\right).
\]
In particular, $\psi_\tau$ preserves the domain $\HH_+$ when $\im(\tau) \ge 0$ and if 
$\alpha\in H_+$, then we have that
$\lim_{\im(\tau)\to\infty}\psi_\tau(\alpha)\cdot\overline{\psi_\tau(\alpha)}=-\infty$.
For the topology on $\widehat H_+$ the convergence of 
$(\alpha_i\in H_+)_i$ to $\alpha_\infty$ means 
that there exist a sequence $(\gamma_i\in \G)_i$ which leave $K$ invariant
and fix $H/K$ pointwise and 
a sequence $(\eta_i\in\RR)_i$ converging to $+\infty$ 
such that $(\alpha_i':=\psi_{\eta_i\sqrt{-1}}^{-1}\gamma_i^{-1}\alpha_i)_i$ converges to 
some $\alpha'_\infty\in H_+$. Then $\alpha_\infty=\pi_K(\alpha'_\infty)$
and
\[
\alpha_i\cdot\bar\alpha_i=\psi_{\eta_i\sqrt{-1}}\alpha'_i\cdot
\overline{\psi_{\eta_i\sqrt{-1}}\alpha'_i}=
\alpha'_i\cdot\bar\alpha'_i -4\eta_i\im 
\left((\alpha'_i\cdot e_0)(\alpha'_i\cdot e_1)\right),
\]
which indeed tends to $-\infty$ as $i\to\infty$.

We next do the case when $K$ has an associated isotropic line $J$:
$K=K_\sigma$ for some $\sigma\in\Sigma_J$. Let $e\in J(\QQ)$ be a
generator. 
For $f\in J^\perp$, we put 
\[
\psi_{e, f}: \alpha\in H\mapsto \alpha +(\alpha \cdot e)f
-(\alpha \cdot f)-\tfrac{1}{2}(f\cdot f)(\alpha\cdot e)e\in H
\]
This transformation respects the quadratic form on $H$ and we have
\[
\psi_{e, f}(\alpha)\cdot\overline{\psi_{e, f}(\alpha)}=
\alpha\cdot\alpha + 
4|\alpha\cdot e|^2 
\left( p_e(\alpha)\cdot\im (f)
+\tfrac{1}{2}\im (f)\cdot\im (f)\right).
\]
So if $f$ is such that $\im (f)\in C_e$, then $\psi_{e, f}$ preserves $H_+$
and $\psi_{e, f}(\alpha)\cdot\overline{\psi_{e, f}(\alpha)}\to -\infty$ if
$\im (f)$ tends to infinity along a ray in $C_e$.
For the topology on $\widehat H_+$ the convergence of 
$(\alpha_i\in H_+)_i$ to $\alpha_\infty\in \pi_KH_+$ means 
that there exist a sequence $(\gamma_i\in \G)_i$ which leave $\sigma$ 
invariant and fix $H/K$ pointwise and 
a sequence $(y_i\in \sigma)_i$ converging to a point $y_\infty$ at infinity
in the relative interior of the projectivized $\sigma$ such that 
$(\alpha_i':=\psi_{e,y_i\sqrt{-1}}^{-1}\gamma_i^{-1}\alpha_i)_i$ converges 
to some $\alpha'_\infty\in H_+$. Then $\alpha_\infty=\pi_K(\alpha'_\infty)$
and
\[
\alpha_i\cdot\bar\alpha_i=\psi_{e,y_i\sqrt{-1}}\alpha'_i\cdot
\overline{\psi_{e,y_i\sqrt{-1}}\alpha'_i}=
\alpha'_i\cdot\bar\alpha'_i +4|\alpha'_i\cdot e|^2 
\left( p_e(\alpha'_i)\cdot y_i +\tfrac{1}{2} y_i\cdot y_i\right)),
\]
which tends to $-\infty$ as $i\to\infty$.
\end{proof}

Let for $\sigma\in\tilde\Sigma$, $\G^\sigma$ stand for the group of $\g\in\G$ 
that leave the stratum $\pi_{K_\sigma}(H_+)$ of $\widehat{H}_+$ pointwise fixed. This is also the group of $\g\in\G$ which fix $K_\sigma^\perp$ pointwise and leave $\sigma$ invariant. In certain cases these subgroups separate the strata:

\begin{lemma}\label{lemma:invar}
If the form  $\alpha\in H\mapsto \alpha\cdot\bar\alpha$ takes on every 
nonisotropic two-dimensional  intersection of members of $\Hcal$ a positive value, then 
for all $\sigma\in\tilde\Sigma$, $K_\sigma^\perp$ is 
the fixed point set of $\G^\sigma$ in $H$. In particular, $\pi_{K_\sigma}(H_+)$ is the fixed point set of $\G^\sigma$ in $\widehat{H}_+$.
\end{lemma}
\begin{proof}
Denote by $G_{K^\perp}$ the group of 
orthogonal transformations of $H$ which leave $K^\perp$ pointwise fixed. 
Suppose first that $K=K_\sigma$ is of type 1, so of signature $(\dim
K-2,2)$. Then by definition
$\G^\sigma$ is the group of $\g\in \G$ which leave
$K^\perp$ pointwise fixed. The latter is an arithmetic group in 
in the orthogonal group of $K$. Since $\dim K\ge 3$, this orthogonal group 
is not anisotropic. This implies \cite{borel} that it contains 
$\G_{K^\perp}$ as a Zariski dense subgroup. 

Next we do the case when $K$ is of type 2. Then again $\G^\sigma$ is 
arithmetic in $G_{K^\perp}$. In particular, it meets the unipotent radical 
$U(G_{K^\perp})$ of 
the latter in an arithmetic subgroup. We observe that $U(G_{K^\perp})$
is the group of unipotent orthogonal transformations of $H$ 
which leave $K^\perp$ pointwise fixed and has $K^\perp$ as its fixed point 
set in $H$. Hence this is also the fixed point set of $\G^\sigma$.

Finally, we do the case when $\sigma\in\Sigma_J$, with $J$ a $\QQ$-isotropic
line. So $K\supset J$ and $K/J$ has hyperbolic signature. 
Choose a generator $e$ of $J(\QQ)$ so that we have defined the distinguished quadratic 
cone $C_e\subset (J^\perp/J)(\RR)$ and $\Sigma_J$ is a locally rational 
polyhedral decomposition of that cone. The span of 
$\sigma$ is $K/J$. Our assumption implies that no member of
$\Sigma_J$ is of dimension one. So the one dimensional rays on the boundary of
$\sigma$ must all be all improper, that is, lie on $\QQ$-isotropic lines
in $C_e$. In particular, the collection $\Lcal$ of such lines spans $K/J$. 
If $L\in \Lcal$, then its preimage in
$\tilde L$ in $H$ is a $\QQ$-isotropic plane which contains $J$. Let $e_L\in
\tilde L(\QQ)$ be such that $(e,e_L)$ is basis of $\tilde L$. Then the
transformation 
\[
\psi_{e,e_L} :H\to H,\quad  \psi_L(\alpha)=\alpha +(\alpha\cdot e)e_L
-(\alpha\cdot e_L)e.
\]
has $\tilde L^\perp$ as its fixed point set. It  
also acts trivially on $J^\perp/J$ and hence does so on $K/J$.
Some power of $\psi_{e,e_L}$ will lie in
$\G$ and so upon replacing $e_L$ be a positive multiple, we may assume that
$\psi_{e,e_L}\in\G$. It is clear that then $\psi_{e,e_L}\in\G^\sigma$. 
Hence the fixed point set of $\G^\sigma$ equals $\cap_{L\in\Lcal} \tilde
L^\perp= (\sum_{L\in\Lcal} \tilde L)^\perp=K^\perp$.
\end{proof}

We abbreviate 
\[
M:=\G\bs \PP(H_+^\circ{}),\quad \widehat{M}:=
\G\bs \PP(\widehat H_+^\circ),\quad M_\infty:=\widehat{M}-M
\]
and denote the inclusion $M\subset \widehat{M}$ by $j$.  The variety 
$\widehat{M}$ is naturally
stratified into orbifolds.

If $\G$ acts neatly on $\PP(H^\circ_+)$ in the sense of Borel, then the strata
of $\widehat{M}$ are smooth and $\Fcal$ is a genuine line bundle.  Since there is always a subgroup of $\G$ of finite index which is neat, we will assume that $\G$ already has that property. 
Then the trivial local system $\tilde\HH:=H_{\PP(H^\circ_+)}$ 
over $\PP(H^\circ_+)$ is in a tautological manner a polarized variation of Hodge structure of weight zero (whose classifying map is the identity). The group $\G$ acts on it so that we get a  variation of Hodge structure 
$\HH$ over $M$ with the line bundle  $j^*\Fcal \subset \Ocal_M\otimes \HH$ defining its Hodge flag of level $1$.

\begin{theorem}\label{thm:coext}
Assigning to  $\sigma\in\tilde\Sigma$ the subspace $(H/K_\sigma)^*$ 
of $H^*$ defines a subsheaf $\VV\subset j_*\HH$ that is a
boundary extension of $\HH$.  The image of 
$\VV\to j_*(\Fcal^1(\HH)^*)$ spans a line bundle whose 
$\Ocal_{\widehat{M}}$-dual can be identified with $\Fcal$. 

If $\dim H\ge 5$ and the hermitian form takes on every two-dimensional intersection of members of $\Hcal$ a positive value, then $M_\infty$ is everywhere of codimension $\ge 2$ in $\widehat M$ and $\VV= j_*\HH$.
\end{theorem}
\begin{proof}
A constructible subsheaf $\tilde \VV$ of the direct image of  
$H^*_{\PP(H^\circ_+)}$ on $\PP(\widehat H_+^\circ)$ is defined by letting it on the stratum $\PP\pi_{K_\sigma}\PP(H_+^\circ)$ 
be constant equal to  $(H/K_\sigma)^*$; this is indeed a sheaf because an incidence relation between strata implies the an  inclusion relation 
of subspaces of $H^*$. Passage to the $\G$-quotient then  defines constructible subsheaf $\VV$ of $j_*\HH$, which is clearly defined over $\QQ$. The assertion concerning 
$\Fcal$ is essentially tautological since the restriction of $\Lcal$ to a stratum associated to $\sigma$ is the image of the $\CC^\times$-bundle $\pi_{K_\sigma}(H_+^\circ)\to\PP\pi_{K_\sigma}\PP(H_+^\circ)$) under formation of the $\G$-quotient.
It remains to see that the norm on $\Fcal^*$ vanishes on
$M_\infty$. This follows from Lemma \ref{lemma:convergence}.

Suppose now that no intersection of members of $\Hcal$ is a plane 
on which the form is negative semidefinite. Clearly, the boundary 
$\widehat M-M$ is everywhere of codimension 
$\ge 2$ precisely when there are no strata of dimension $n-1$. 
This is the case if and only if no member of $\tilde\Sigma$ 
is of dimension $2$, or equivalently, if $\dim H\ge 5$ and no 
nonzero intersection of members of $\Hcal$ is negative semidefinite.
The equality $\VV= j_*\HH$ follows from Lemma \ref{lemma:invar}.
\end{proof}

The preceding was also carried out in a ball quotient setting.
Things then simplify considerably. This case may arise as a restriction of the case discussed above as follows.
Suppose that $H$ comes with a faithful action of the group $\mu_l$ of $l$-th 
roots of unity with $l\ge 3$, such
that if $\chi:\mu_l\subset\CC^\times$ is the tautological character, 
the hermitian form  has on the $\chi$-eigen space $H_\chi$ hyperbolic signature (i.e., $(1,m)$ for some $m\ge 0$), or equivalently, that the orthogonal complement of the 
$\RR$-subspace $H_\chi +H_{\bar\chi}$  is positive definite.

We observed that $H_\chi$  is isotropic for the quadratic form so that
intersection $H_{\chi,+}:=H_\chi\cap H_+$ is simply 
the open set of $\alpha\in H_\chi$ with $\alpha \cdot\bar\alpha<0$. 
This implies that $\PP(H_{\chi,+})$ is a complex 
$m$-ball. Now let $\G'$ be an arithmetic subgroup in the 
centralizer of $\mu_l$ in the orthogonal group of $H$. 
Then $\G'$ acts properly on $\PP(H_{\chi,+})$ and with finite covolume.

Let $\Hcal'$ be a collection of $\mu_l$-invariant hyperplanes 
of $H_\chi$ of hyperbolic signature that are obtained by intersecting $H_\chi$ with
a hyperplane in $H$ defined over $\QQ$. We also assume that $\Hcal'$
is a union of a finite number of $\G'$-equivalence classes. The open subset   
$H_{\chi,+}^\circ\subset H_{\chi,+}$ is defined as usual, but the role of 
$\tilde\Sigma$ is now taken by a selfindexed collection
$\Kcal'=\Kcal'_{1}\cup\Kcal'_{2}$ (the fact that there are no subspaces 
of type 3 is the reason that things simplify). 
Here $\Kcal'_{1}$ is the collection
of intersections of members of $\Hcal'$ of hyperbolic signature and
$\Kcal'_{2}$ is bijectively labeled by the collection of the isotropic 
lines $J'\subset H_{\chi}$ for which $J=J'+\bar J'$ is defined over $\QQ$ (this
is then  a $\QQ$-isotropic plane): for each such $J'$ we let $K_{J'}$ be the 
intersection of $J^\perp\cap H_\chi$ with the collection of members of $\Hcal'$
which contain $J'$. From this point onwards one proceeds as in the case
considered above. We obtain a projective compactification of 
$M':=\G'\bs H_{\chi,+}$, $M'\subset \widehat M'$ and an ample extension $\Fcal'$ over 
$\widehat{M}'$ of the automorphic line bundle over $M'$ and find:

\begin{theorem}\label{thm:coext2'}
There is a canonical boundary extension of Hodge type
$\VV'\subset j'_*\HH_\chi$ such that the image of 
$\VV'\to j'_*(\Fcal^1(\HH)^*)$ spans a line bundle whose 
$\Ocal_{\widehat M'}$-dual can be identified with $\Fcal'$. 

If $\dim H_\chi\ge 3$ and every one-dimensional intersection 
of members of $\Hcal'$  is positive for the hermitian form, then the boundary 
$M'_\infty:=\widehat M'-M'$ is everywhere of codimension $\ge 2$ in 
$\widehat M'$ and $\VV'= j'_*\HH_\chi$.
\end{theorem}

\end{document}